\documentclass[12pt]{amsart}

\usepackage{amsmath,amssymb,amsthm,enumerate}

\newtheorem{theorem}{Theorem}[section]

\newtheorem{proposition}[theorem]{Proposition}

\theoremstyle{definition}

\theoremstyle{remark}
\newtheorem{remark}[theorem]{Remark}

\newcommand{\bD}{\mathbb D}
\newcommand{\bP}{\mathbb P}
\newcommand{\bC}{\mathbb C}
\newcommand{\bR}{\mathbb R}

\newcommand{\bZ}{\mathbb Z}

\newcommand{\beq}{\begin{equation}}
\newcommand{\eeq}{\end{equation}}

\newcommand{\eps}{\epsilon}

\DeclareMathOperator{\im}{Im}

\begin{document}
\title[]{A family of compact strictly pseudoconvex hypersurfaces in $\bC^2$ without umbilical points}
\author{Peter Ebenfelt}
\address{Department of Mathematics, University of California at San Diego, La Jolla, CA 92093-0112}
\email{pebenfel@math.ucsd.edu}
\author{Duong Ngoc Son}
\address{Texas A\&M University at Qatar, Science Program, PO Box 23874, Education City, Doha, Qatar}
\email{son.duong@qatar.tamu.edu}

\author{Dmitri Zaitsev}
\address{Department of Mathematics, Trinity College, Dublin}
\email{zaitsev@maths.tcd.ie}
\date{\today}
\thanks{The first author was supported in part by the NSF grant DMS-1600701.}
\thanks{The second author was partially supported by the Qatar National Research Fund, NPRP project 7-511-1-098. Part of this work was done while the second author visited UC San Diego, which he would like to thank for partial support and hospitality.}
\begin{abstract}
We prove the following: {\it
For $\eps>0$, let $D_\eps$ be the bounded strictly pseudoconvex domain in $\bC^2$ given by
\begin{equation*}
(\log|z|)^2+(\log|w|)^2<\eps^2.
\end{equation*}
 The boundary $M_\eps:=\partial D_\eps\subset \bC^2$ is a compact strictly pseudoconvex CR manifold without umbilical points.} This resolves a long-standing question in complex analysis that goes back to the work of S.-S. Chern and J. K. Moser in 1974.
\end{abstract}


\maketitle

\section{Introduction}

A long-standing and well known question concerning the geometry of domains and their boundaries in several complex variables is the following:
\smallskip

{\it Does there exist a bounded strictly pseudoconvex domain $D\subset\bC^2$ with smooth boundary $M:=\partial D$ such that $M$ has no (CR) umbilical points?}
\smallskip

This question, which originated in the seminal 1974 paper \cite{CM74} by S.-S. Chern and J. K. Moser and is often referred to as the Chern--Moser question, is particular to $\bC^2$. It became clear already in \cite{CM74} that there is sharp difference between umbilical points on strictly pseudoconvex boundaries in $\bC^2$ and in $\bC^n$ with $n\geq 3$. To begin with, umbilical points in $\bC^n$ with $n\geq 3$ are determined by the vanshing of a 4th order tensor, whereas umbilical points in $\bC^2$ are determined by the vanishing of a 6th order tensor (discovered already by E. Cartan \cite{Cartan32,Cartan33} in the early 1930's). A simple Thom transversality argument (see, e.g., \cite{EZ16}) shows that a generic (i.e., sufficiently general) strictly pseudoconvex domain in $\bC^n$ with $n\geq 4$ does not have any umbilical points in its boundary, and Webster \cite{Webster00} showed that, in particular, every non-spherical real ellipsoid in $\bC^n$ with $n\geq 3$ has no umbilical points. In contrast with Webster's result, and illustrating the point that the situation in $\bC^2$ and that in $\bC^n$ with $n\geq 3$ is different, X. Huang and S. Ji \cite{HuangJi07} proved that every real ellipsoid in $\bC^2$ must have umbilical points. We mention here also two other recent papers, \cite{EDumb15} and \cite{EZ16}, in which the focus (regarding Chern--Moser's question) has been on proving results to the effect that certain classes of three-dimensional CR manifolds must possess umbilical points (supporting a possible 'no' as an answer to Chern--Moser's question). In the former, it is shown, e.g., that boundaries of bounded, complete circular domains must have umbilical points. In the latter, it is shown, e.g., that generic "almost circular" perturbations of the unit sphere also must possess umbilical points.

The purpose of this note is to settle Chern--Moser's question by proving the following:

\begin{theorem}\label{T:main0}
For $\eps>0$, let $D_\eps$ be the bounded strictly pseudoconvex domain in $\bC^2$ given by
\beq\label{Deps}
(\log|z|)^2+(\log|w|)^2<\eps^2.
\eeq
 The boundary $M_\eps:=\partial D_\eps\subset \bC^2$ is a compact strictly pseudoconvex CR manifold without umbilical points.
\end{theorem}

\begin{remark} We make here two remarks:

\begin{itemize}
\item[(a)]
The domains $D_\eps$ are biholomorphically inequivalent for different $\eps>0$, and the boundaries $M_\eps$ are not homogeneous. However, the boundaries are locally homogeneous and locally equivalent. This is explained in more detail in Section \ref{S:compact}.
\item[(b)] We note that the $D_\eps$ are Reinhardt domains, and hence circular. However, they are not {\it complete} circular, and no circle action on $D_\eps$ is everywhere transversal (to the CR structure). Thus, Theorem \ref{T:main0} does not contradict the existence results in \cite{EDumb15}.
\end{itemize}
\end{remark}

The motivation behind this example comes from the context of Grauert tubes of Riemannian manifolds. This context is only conceptual and will play no role in the details given in this paper. We will therefore refrain from giving any formal definition or general discussion of this notion; readers with a knowledge of Grauert tubes would do well to recognize and keep this context in mind, and the interested reader without such knowledge may consult \cite{LempertSzoke91} or \cite{GuilleminStenzel91}. We first recall that in
his classification of homogeneous three-dimensional strictly pseudoconvex CR manifolds, E. Cartan \cite{Cartan32,Cartan33} discovered that the homogeneous, compact and strictly pseudoconvex CR manifolds
\begin{equation}\label{mualph}
\mu_\alpha:=\{[z_0:z_1:z_2]\in \bC\bP^2\colon |z_0|^2+|z_1|^2+|z_2|^2=\alpha
|z_0^2+z_1^2+z_2^2|\},\quad \alpha >1,
\end{equation}
and their covers, which were subsequently classified in \cite{Isaev06} as a $4\!:\! 1$ cover $\mu_\alpha^{(4)}$ (diffeomorphic to a sphere) that factors through a $2\!:\! 1$ cover $\mu_\alpha^{(2)}$ (consisting of the intersection of sphere and a nonsingular holomorphic quadric in $\bC^3$), all have no umbilical points, and they are the only possible homogeneous and compact examples. It is known, however, that none of these embeds in $\bC^2$ (see below). The first observation (see \cite{PatrizioWong91}) is then that Cartan's family of homogeneous examples $\mu_\alpha$ occur as boundaries of Grauert tubes of the sphere $S^2$ with the round (standard) metric, and the $2\!:\! 1$ covers $\mu_\alpha^{(2)}$ occur as the boundaries of Grauert tubes of the real projective plane $\bR\bP^2$ with its constant curvature metric. These tubes cannot be biholomorphically embedded as domains in $\bC^2$, because if they could, then $S^2$ and $\bR\bP^2$ would be embedded as totally real submanifolds in $\bC^2$, which is not possible; see \cite{Bishop65}, \cite{Wells69}, \cite{NemirovskiSiegel16}. (The $4\!:\! 1$ cover $\mu_\alpha^{(4)}$ is not a Grauert tube and cannot be embedded in $\bC^n$ for any $n$; indeed, the $\mu_\alpha^{(4)}$ were given by Rossi \cite{Rossi65} as examples of CR manifolds that cannot be so embedded; see also \cite{Isaev06}.) With this perspective, the next potential example that comes to mind then is a Grauert tube over the flat torus or the flat Klein bottle. These can indeed be biholomorphically embedded as bounded domains in $\bC^2$ (at least for small radii) by complexifying any real-analytic, totally real embedding of the torus or Klein bottle. It turns out that, in fact, the Grauert tubes over the flat torus or Klein bottle have boundaries without umbilical points (see Section \ref{S:compact}). The family of examples in Theorem \ref{T:main0} is the family of Grauert tubes over a flat torus, embedded as domains in $\bC^2$, and consequently have the properties described in the theorem.

This note is organized as follows. In Section \ref{S:noncompact}, a homogeneous non-compact model without umbilical points is discussed. In Section \ref{S:compact}, the Grauert tubes over the torus and Klein bottle are constructed from the noncompact model, and an embedding of the Grauert tube over the torus is explicitly given, yielding the family of examples in Theorem \ref{T:main0}. In Section \ref{S:concluding}, a possible refinement of Chern--Moser's question is proposed.

\subsection{Acknowledgements} The authors would like to thank Stefan Nemirovski and Alexander Sukhov for useful discussions and the suggestion to consider domains with non-trivial topology.

\section{The noncompact model}\label{S:noncompact}

Consider the strictly plurisubharmonic polynomial $\rho$ in $\bC^2$ given by
\beq\label{rho}
\rho=\rho(z,w,\bar z,\bar w):=(\im z)^2+(\im w)^2.
\eeq
For $\epsilon>0$, we shall consider the strictly pseudoconvex domain $\tilde\Omega_\eps\subset \bC^2$ given by
\beq\label{Omega}
\tilde \Omega_\eps:=\{(z,w)\colon \rho(z,w,\bar z,\bar w)<\eps^2\}
\eeq
and its boundary
\beq\label{Omega}
\tilde M_\eps:=\{(z,w)\colon \rho(z,w,\bar z,\bar w)=\eps^2\}.
\eeq
A direct calculation shows that $i\partial\bar \partial\rho=i(dz\wedge d\bar z +dw\wedge d\bar w)/2$, and $(\partial \bar\partial \sqrt{\rho})^2=0$, which means (although this will not be important here) that $\Omega_\eps$ is the Grauert tube of radius $\eps$ of the flat Euclidian metric on $\bR^2$. These hypersurfaces were discovered by E. Cartan \cite{Cartan32,Cartan33} in his classification of homogeneous, strictly pseudoconvex three-dimensional CR manifolds. For each $\eps>0$, the automorphism group of $\tilde M_\eps$ is given by $O(2,\bR)\times \bR^2$ via the action
$$
(z,w)\mapsto (z,w)A+(a,b),\quad A\in O(2,\bR),\ (a,b)\in \bR^2.
$$
Moreover, a simple scaling $(z,w)\mapsto (rz,rw)$ shows that any two $\tilde M_\eps$, $\tilde M_{\eps'}$ are biholomorphically equivalent. Cartan also showed that these hypersurfaces are {\it non-spherical}, which as we recall for a {\it homogeneous} CR manifold is equivalent to having no umbilical points. We summarize these observations in the following:

\begin{proposition}\label{P:tildenonumbilic}
The strictly pseudoconvex hypersurface $\tilde M_\eps$, for any $\eps>0$, is homogeneous, strictly pseudoconvex, and non-spherical (i.e., has no umbilical points). Moreover, $\tilde M_\eps$ and $\tilde M_{\eps'}$ are biholomorphically equivalent for any pair $\eps,\eps'>0$.
\end{proposition}

As mentioned above, the proof of Proposition \ref{P:tildenonumbilic} is contained in the work of Cartan, {\it op.cit.}. For the reader's convenience, however, we give here a proof of the key fact that $\tilde M_\eps$ has no umbilical points. To this end, we shall utilize a recent new characterization of umbilical points due to the first and third authors in  \cite{EZ16}: it suffices (see, e.g., formula (3.5), {\it Ibid.}) to check that
\beq\label{deta3}
\det A_3(\rho):=
\det \begin{pmatrix}
\rho_w^3  & \bar L(\rho_w^3) & \cdots & \bar L^4(\rho_w^3)\cr
 \rho_z\rho_w^2 & \bar L( \rho_z\rho_w^2) & \cdots & \bar L^4( \rho_z\rho_w^2)\cr
  \rho_z^2\rho_w & \bar L(  \rho_z^2\rho_w) & \cdots & \bar L^4(  \rho_z^2\rho_w)\cr
\rho_z^3 & \bar L(\rho_z^3) & \cdots & \bar L^4(\rho_z^3)\cr
\rho_{Z^2}(L, L) & \bar L(\rho_{Z^2}(L, L)) & \cdots & \bar L^4(\rho_{Z^2}(L, L))
\end{pmatrix}
\eeq
does not vanish on $\tilde M_\eps$. Here $L$ is the (1,0) vector field
\beq
L:=-\rho_w\frac{\partial}{\partial z}+\rho_z\frac{\partial}{\partial w},
\eeq
which is tangent to $\tilde M_\eps$, for every $\eps>0$. We note that
\beq
\rho_z=-\frac{1}{2}(z-\bar z),\quad \rho_w=-\frac{1}{2}(w-\bar w), \quad \rho_{Z^2}(L,L)=-\frac{1}{8}((z-\bar z)^2+(w-\bar w)^2)=\frac{1}{2}\rho.
\eeq
Since $\bar L$ annihilates $\rho$, the only non-zero entry in the last row of $A_3(\rho)$ is $\rho_{Z^2}(L,L)=\rho/2$, and we conclude that on $\tilde M_\eps$
\beq\label{deta3-2}
\det A_3(\rho)|_{\rho=\eps^2}=\frac{1}{2}\eps^2\det B|_{\rho=\eps^2},
\eeq
where $B$ is the matrix
\beq\label{Bmatrix}
B:=
\begin{pmatrix}
\bar L(\rho_w^3) & \cdots & \bar L^4(\rho_w^3)\cr
  \bar L( \rho_z\rho_w^2) & \cdots & \bar L^4( \rho_z\rho_w^2)\cr
   \bar L(  \rho_z^2\rho_w) & \cdots & \bar L^4(  \rho_z^2\rho_w)\cr
 \bar L(\rho_z^3) & \cdots & \bar L^4(\rho_z^3)
\end{pmatrix}.
\eeq
We observe that the following identities hold:
\beq\label{rhoder}
\frac{\partial}{\partial \bar z}\, \rho_z=\frac{1}{2}, \quad \frac{\partial}{\partial \bar w}\,\rho_z=0,\quad \frac{\partial}{\partial \bar z}\,\rho_w=0,\quad \frac{\partial}{\partial \bar w}\,\rho_w=\frac{1}{2}.
\eeq
We also note, in view of the fact that each $\tilde M_\eps$ is homogeneous and any two $\tilde M_\eps$, $\tilde M_{\eps'}$ are biholomorphically equivalent, that it suffices to verify that $\det B=\det B(z,w,\bar z,\bar w)$ is non-zero at {\it one single point}. If we choose, e.g., the point $(z,w)=(0,1)$ and use \eqref{deta3-2}, \eqref{rhoder}, then it is easily verified (and left to the reader) that $\det B|_{(0,1)}\neq 0$. Thus, we may conclude that
$\det A_3(\rho)$ does not vanish outside $\bR^2=\{\rho=0\}\subset \bC^2$. (One may in fact show, although we shall not do so here, that $\det A_3(\rho)$ equals a non-zero constant times $\eps^{14}$ on each $\tilde M_\eps$ by interpreting the matrix $B$ as a Wronskian.) This completes the proof that $\tilde M_\eps$ has no umbilical points for any $\eps>0$.

\section{The compact tubes over flat $\bR^2$-quotients}\label{S:compact}

Let $\Lambda$ be a discrete subgroup in the group of rigid motions of $\bR^2$ such that $\Sigma=\Sigma^2:=\bR^2/\Lambda$ is a compact (smooth) 2-surface. For example, if $\Lambda$ is generated by the two translations $(x,y)\mapsto (x,y)+e_j$, $j=1,2$, where $e_1,e_2$ are linearly independent vectors in $\bR^2$, then $\Sigma$ is a 2-torus. If $\Lambda$ is generated by, e.g.,  $(x,y)\mapsto (x,y+1)$, $(x,y)\mapsto (x+1,-y)$, then $\Sigma$ is the Klein bottle. We observe that the function $\rho$ in \eqref{rho} is well defined (descends) on the complex manifold $\Omega:=\bC^2/\Lambda$, and that $\Sigma$ sits as a totally real 2-surface in $\Omega$. (This observation goes back to \cite{LempertSzoke91}.) We shall consider the strictly pseudoconvex tube domains
\beq
\Omega_\eps:=\{(z,w)\in \Omega\colon \rho(z,w,\bar z,\bar w)<\eps^2\},
\eeq
for $\eps>0$, and their boundaries $M_\eps:=\partial\Omega_\eps$. We note, although this observation will not be important here, that $\Omega_\eps$ is the Grauert tube of radius $\eps>0$ over the 2-surface $\Sigma$ equipped with the flat metric (see, e.g., \cite{LempertSzoke91}, Example 2.1). The CR manifold $M_\eps$ is compact and, clearly, locally biholomorphic to the noncompact CR manifold $\tilde M_\eps$ in the previous section. Thus, it follows from Proposition \ref{P:tildenonumbilic} that $M_\eps$ is a compact, strictly pseudoconvex CR manifold without umbilical points.

We mention that the compact CR manifolds $M_\eps$ obtained in this way are not homogeneous, in contrast with their non-compact models $\tilde M_\eps$. As a consequence of a rigidity result of Burns and Hind for Grauert tubes over compact manifolds \cite{BurnsHind01}, the automorphism group of $M_\eps$ is isomorphic, via the tangent map, to the group of isometries of the 2-surface $\Sigma$, which, e.g., in the flat torus case is 2-dimensional. Moreover, it follows from the same rigidity result that the CR manifolds $M_\eps$, $M_{\eps'}$ for $\eps\neq \eps'$ are not equivalent. However, since $\tilde\Omega_\eps$ is the universal cover of $\Omega_\eps$ and the covering map extends analytically to the boundary $\tilde M_\eps\to M_\eps$, we may conclude that for any $\eps,\eps'>0$ and $p\in M_\eps$, $p'\in M_{\eps'}$, there is a local biholomorphism such that, locally, $(M_\eps,p)\cong (M_{\eps'},p')$. Thus, we have {\it local}, but not global, homogeneity and equivalence for the family of compact CR manifolds $M_\eps$.

We now note that if $\Sigma$ has a real-analytic, totally real embedding $f\colon \Sigma\to \bC^2$ (e.g., when $\Sigma$ is either a 2-torus or the Klein bottle \cite{Rudin81}), then by compactness and complexification of the real-analytic map $f$, we deduce that there exists a sufficiently small $\eps_0>0$ and a holomorphic map $F\colon \Omega_{\eps_0}\to\bC^2$ such that $F|_\Sigma=f$. Since $f$ is a totally real embedding, we also conclude that for sufficiently small $\eps_0$, the holomorphic map $F$ is a biholomorphism onto its image $F(\Omega_{\eps_0})\subset \bC^2$. The following is then a direct consequence of this construction:

\begin{theorem}\label{T:main1}
With notation as above, the real hypersurfaces $F(M_\eps)\subset \bC^2$, for $0<\eps<\eps_0$, are compact, strictly pseudoconvex hypersurfaces (each bounding the domain $F(\Omega_\eps)\subset \bC^2$) that do not have any umbilical points.
\end{theorem}

We also note that, from the Grauert tube construction, it follows that each manifold $M_{\eps}$ is a circle bundle in the tangent bundle of the 2-surface $\Sigma$ with flat
metric and the CR structure is induced by the \emph{adapted complex structure}. This class of CR manifolds is different from the unit circle bundles in hermitian
line bundles considered in \cite{EDumb15}; manifolds from latter class, in the case of the 2-torus $\Sigma$, with two-dimensional automorphism groups must have umbilical points, by \cite[Theorem~1.4]{EDumb15}.

\subsection{An explicit construction for the Grauert tube over the standard 2-torus. Proof of Theorem $\ref{T:main0}$.} We shall now give an explicit construction of the embedded Grauert tubes $F(\Omega_\eps)$ over the standard 2-torus $\Sigma=\mathbb T^2:=\bR^2/2\pi\bZ^2$; thus, the group $\Lambda$ is generated by $(x,y)\mapsto (x+2\pi,y)$ and $(x,y)\mapsto (x,y+2\pi)$ and  $\rho$ (given by \eqref{rho}) corresponds to the flat Euclidian metric on $\mathbb T^2$. In this case, we have an explicit global totally real embedding $f\colon \Sigma\to \bC^2$ given by
\beq\label{TRemb}
f(x,y):=(e^{ix},e^{iy}),
\eeq
which complexifies to a global holomorphic embedding $F\colon \bC^2/\Lambda\to \bC^2$:
\beq\label{CTRemb}
F(X,Y):=(e^{iX},e^{iY}),\quad X=x+ix',\, Y=y+iy'.
\eeq
In the coordinates
\beq
z=e^{iX},\quad w=e^{iY},
\eeq in $\bC^2$, the embedded Grauert tube $D_\eps:=F(\Omega_\eps)$ is then given by the equation:
\beq
(\log|z|)^2+(\log|w|)^2<\eps^2.
\eeq
Theorem \ref{T:main0} is now a direct consequence of Theorem \ref{T:main1}. \qed


\section{Concluding remarks}\label{S:concluding}

We observe that the Grauert tube over the flat torus (i.e., the example $D_\eps$ in Theorem \ref{T:main0}) is a disk bundle in the tangent bundle over the torus, which is a trivial bundle. Hence, $D_\eps$ is diffeomorphic to $\mathbb T^2\times \bD$, where $\bD$ denotes the unit disk in $\bC$. One could refine the Chern--Moser question in the introduction, so as to make it open again, by asking if the domain $D\subset \bC^2$ could be found such that {\it $\overline D$ is also diffeomorphic to the closed unit ball in $\bC^2$} (making the boundary diffeomorphic to the 3-sphere $S^3$). We remark that the only known (to the authors) examples of compact three-dimensional strictly pseudoconvex CR manifolds without umbilical points that are diffeomorphic to $S^3$ are the Rossi spheres $\mu^{(4)}_\alpha$ and their perturbations, and these do not even have Stein neighborhoods.

We recall that the examples given in this paper are such that their boundaries are locally homogeneous. We remark that if a bounded strictly pseudoconvex domain $D\subset \bC^2$ is such that $\overline D$ is diffeomorphic to the closed unit ball in $\bC^2$ and its boundary $M:=\partial D$ has no umbilical points, then $M$ cannot be locally homogeneous. Indeed, we note that $M$ is connected and simply connected. Hence, it follows from a classical result of Pinchuk \cite{Pinchuk75} that any local automorphism of $M$ can be extended to a global one. Combining this observation with a standard connectedness argument, we conclude that if $M$ is locally homogeneous, it is in fact also (globally) homogeneous. Now, by the classification results of Cartan \cite{Cartan32,Cartan33} and Isaev \cite{Isaev06}, $M$ must then be a Rossi sphere $\mu^{(4)}_\alpha$, and as noted above these cannot be embedded into $\bC^2$.


\def\cprime{$'$}

\end{document}